\documentclass[12pt]{article}
\usepackage{graphicx}
\usepackage{amssymb}
\usepackage{amsmath}
\usepackage{amsthm}
\usepackage{mathtools}
\usepackage{color}
\numberwithin{equation}{section}
\mathtoolsset{showonlyrefs=false}

\newcommand{\R}{{\mathbb R}}

\newtheorem{theorem}{Theorem}

\newtheorem{prop}[theorem]{Proposition}

\newtheorem{oss}{Remark}
\newtheorem{ex}[oss]{Example}
\newtheorem*{theorem*}{Theorem}
\def\e{\varepsilon}

\title{\bf Invariant cones for linear elliptic systems\\ with gradient coupling}
\author{I. Capuzzo Dolcetta$^{\hbox{\small{ a}}}$, L. Rossi$^{\hbox{\small{ a,b }}}$ and  A. Vitolo$^{\hbox{\small{ c}}}$\\
	\\
\footnotesize{$^{\hbox{a }}$Dipartimento di Matematica, Sapienza Universit\`a  di Roma, Italy}\\
\footnotesize{$^{\hbox{b }}$CAMS, CNRS-EHESS, Paris, France}\\
\footnotesize{$^{\hbox{c }}$Dipartimento di  Ingegneria Civile, Universit\`a di Salerno, Italy}}
\date{}

\begin{document}
\maketitle

\bigskip
\bigskip
\bigskip
\noindent {\bf Abstract.} {We discuss counterexamples to the validity of the weak Maximum Principle for linear elliptic systems with zero and first order couplings and prove, through a suitable reduction to a nonlinear scalar equation, a quite general result showing that  some algebraic condition on the structure of gradient couplings and a 
	cooperativity condition on the matrix of zero order couplings guarantee the existence of invariant cones 
	in the sense of Weinberger~\cite{Wein}}.

\noindent {\small {\em MSC 2010 Numbers:} 35J47, 35J70, 35B50, 35P30, 35D40}
\footnote
{This work has been partially supported by GNAMPA-INdAM\\

\noindent i.capuzzodolcetta@gmail.com\\
l.rossi@uniroma1.it\\
vitolo@unisa.it}

\noindent {\small {\em Keywords and phrases: elliptic differential inequalities, weak Maximum Principle, invariant cones, Bellman operators, principal eigenvalue}. 

\section{Introduction}\label{intro}

We consider smooth vector-valued functions $u=(u_1,\dots,u_m)$ of the variable $x$ in a bounded open subset $\Omega\subset \R^n$ satisfying linear systems of partial differential inequalities  of the following form
\begin{equation}\label{systemcomp}
Au +\sum_{i=1}^n B^{(i)} D_iu+Cu \geq 0\;\; \hbox{in}\;\; \Omega
\end{equation}
where $A$ is the second order operator 
\begin{equation}\label{A}
Au=\left(\begin{array}{c}\Delta u_1 \\. \\.\\. \\ \Delta u_m\end{array}\right)
\end{equation}
$B^{(i)}$ and $C$ are $m\times m$ real matrices and with constant coefficients, and for $i=1,\dots,n$,
\begin{equation}\label{D_i}
D_i u=\left(\begin{array}{c}  \frac{\partial u_1}{\partial x_i}
\\. \\.\\. \\ \frac{\partial u_m}{\partial x_i}
\end{array}\right)
\end{equation}
denotes the $i-th$ column of the Jacobian matrix of the vector function $u$.\\ Note that the above defined structure of the systems allows coupling between the $u_j$ and their gradients but not at the level of second derivatives.\\
Specific assumptions on the $B^{(i)}$ and $C$  will be made later on.\\ 

Systems of this kind naturally arise in several different contexts such as modeling of simultaneous diffusions of $m$ substances which decay spontaneously or in the case of systems describing switching diffusion processes in probability theory. In the latter case the homogeneous Dirichlet problem for system (\ref{systemcomp}) describes discounted exit  times from $\Omega$, see for example \cite{ChenZhao}.
%Sirakov-Souplet etc.....\\

We are interested here in investigating the validity of the weak Maximum Principle, {\bf wMP} in short, that is  the sign propagation property from the boundary to the interior for solutions $u=(u_1,\dots,u_m)$ of the differential inequalities (\ref{systemcomp}), i.e.,
\begin{equation}\label{sp}
\mathbf{wMP}: \quad  u\leq 0  \ \mbox{on} \ \partial \Omega \ \ \implies \ \ u\leq 0 \ \mbox{in}\ \Omega.
\end{equation}
The vector function $u$ will be always assumed to belong to $[C^2(\Omega)]^m\cap [C^0(\overline\Omega)]^m$ and we will adopt the standard notation $u\leq 0$ if $u_j\leq 0$ for each $j=1,\dots, m$.
We adopt the same notation for real-valued matrices, namely for a matrix $A$, $A\geq0$ means that all its entries 
are nonnegative.
\\
The validity of {\bf wMP} is well-understood in the scalar case $m=1$ even for 
general degenerate elliptic fully nonlinear partial differential inequalities such as
$$ F(x,u,\nabla u, \nabla^2 u) \geq 0$$
in a bounded $\Omega$ and also in some unbounded domain of $\R^n$, see \cite{CDV1,CDV2} for recent results in this direction. Let us point out that the {\bf wMP} property in the scalar case is related, and in fact equivalent, to the positivity of the principal eigenvalue (may be a pseudo one, if degeneracy occur in the dependence of $F$ with respect to Hessian matrix $\nabla^2 u$) of the Dirichlet problem for $F$,
see~\cite{BCPR},\cite{BNV}.
 
The case $m>1$ has been the object of several papers mainly in the  case of diagonal weakly coupled systems, that is when the matrices $B^{(i)}$ are diagonal and couplings between the functions~$u_j$ only occur at the level of zero-order terms, 
described by a  matrix $C=(c_{jk})_{j,k}$, satisfying the cooperativity condition
\begin{equation} \label{coop} 
c_{jk} \geq 0 \ \ \mbox{for } j\neq k\;, \qquad  \sum_{k=1}^{m}c_{jk} \leq 0 \ \ \mbox{for }\,j=1,\dots, m
\end{equation}
Referring to the aforementioned exit time model, condition (\ref{coop}) requires the discount factor for the $j$-th process to dominate the sum of the interactions coefficients with all the other processes.\\

In the framework of purely weak cooperative couplings, let us mention the results in Section~8 of the book by Protter and Weinberger \cite{PW} and the references therein. For generalizations of those results in some semilinear cases 
see \cite{Sweers},\cite{BMS},\cite{Amann}, while \cite{BS} contains results in the same direction concerning fully nonlinear uniformly elliptic operators $F=F(x,u,\nabla u, \nabla^2 u)$.\\
The recent paper \cite{CDV3} extends the validity  of some of the results in \cite{BS} concerning {\bf wMP}  to a large class of fully nonlinear degenerate elliptic operators.\\ In particular, for the case of linear systems as (\ref{systemcomp}) with no coupling in first derivatives (i.e. when each $B^{(i)}$ is diagonal), 
the main result in  \cite{CDV3}  is that  {\bf wMP} holds true for system (\ref{systemcomp}) provided $C$ is cooperative.\\ 
Let us also point out that the main result in \cite{CDV3} holds even in the more general case where the Laplace operator is replaced by more general expressions $\mathrm{Tr}(A^j \nabla^2 u_j)$ satisfying $\mathrm{Tr}(A^j)>0 $ for $j=1,\dots,m$.\\

When coupling in first order terms occurs in (\ref{systemcomp}), simple examples as the following one taken 
from \cite{BS} show that the {\bf wMP} property (\ref{sp}) may indeed fail: 
\begin{ex}\label{exBS}{\rm The vector $u(x_1,x_2)=(1-x_1^2-x_2^2, \frac{1}{3} x_1^3+4x_2-20)$ is a solution of
	\[
	\begin{cases}
	\displaystyle
	\Delta u_1+ \frac{\partial u_2}{\partial x_2}=0\\
		\displaystyle	
			\Delta u_2+\frac{\partial u_1}{\partial x_1}=0
	
	\end{cases}
	\]
	in the unit ball $\Omega \subset \R^2$, $u_1=0$, $u_2<0$ on $\partial \Omega$ but $u_1>0$ in $\Omega$. Observe that 
the zero-order matrix is $C\equiv 0$ in this example, so that (\ref{coop}) is fulfilled.
}\end{ex}

As a matter of fact, even a first-order coupling of arbitrarily small size in the system
can be responsible of the loss of {\bf wMP}, as the following example shows:

\begin{ex}\label{ex-instability}{\rm The system 
		\[
		\begin{cases}
		\displaystyle
		\Delta u-\varepsilon \frac{\partial v}{\partial x_1} \geq 0\\
		\displaystyle
		\Delta v-\varepsilon' \frac{\partial u}{\partial x_1} \geq 0
		\end{cases}
	\]
		in a bounded domain $\Omega\subset\R^n$, fulfills {\bf wMP} if and only if $\varepsilon=\varepsilon'=0$.\\
		Indeed the validity of {\bf wMP} when $\varepsilon=\varepsilon'=0$ is classical.
		Conversely, if, say, $\varepsilon'\neq0$, then {\bf wMP} is violated by the pair
		$$u(x)=\delta- |x-\bar x|^2,\qquad
		v(x)=v(x_1)=C(e^{-x_1}-H),$$
		where $\bar x\in\Omega$ and $\delta>0$ is small enough to have $u<0$ on $\partial\Omega$,
		and $C\gg1$, $H\gg1$.
}\end{ex}

Example \ref{ex-instability} enlightens an instability property of {\bf wMP} for cooperative systems
with respect to first order perturbations. This is in striking contrast with the scalar case. 
Indeed, for a uniformly elliptic scalar inequality, not only
the presence of a first order term does not affect the validity of {\bf wMP} when the zero-order term
is nonpositive, but in addition {\bf wMP} is stable with respect to perturbations of the coefficients, 
in the $L^\infty$ norm. This can be seen as a consequence of the fact that {\bf wMP}
is characterized by the positivity of the associated principal eigenvalue, and the latter 
depends continuously on the coefficients of the operator, see e.g.~\cite{PW,BNV} and also
\cite{BMS,BS} where such characterization in terms of the same notion of principal 
eigenvalue as in~\cite{BNV} is extended to cooperative systems without first-order coupling.
Example \ref{ex-instability} reveals either that such notion does not exist 
when there is a first-order coupling, or that it is not continuous with respect to the coefficients.

According to the above considerations, two perspectives can be adopted 
in order to investigate the sign-propagation properties for coupled systems 
such as~\eqref{systemcomp}. The first one consists in strengthening the hypotheses on the coefficients 
of the operators, namely the cooperativity condition~\eqref{coop}. 
The second one is to replace {\bf wMP} by some different kind of propagation property which reflects 
in some way the geometry of the coupling terms. We will explore both directions. 

Observe that the systems in Examples~\ref{exBS} and \ref{ex-instability} fulfill the cooperativity condition~\eqref{coop} 
in the ``border case'', that is, when all inequalities are replaced by equalities.
A natural question is then whether it is possible, for the {\bf wMP} to hold, to allow some 
coupling in first-order terms  in the system, at least when the cooperativity conditions~\eqref{coop} 
hold with strict inequalities, i.e.,
\begin{equation}\label{coop>}
	c_{jk} \geq K \ \ \mbox{for } j\neq k\;, \qquad  \sum_{k=1}^{m}c_{jk} \leq -K \ \ \mbox{for }\,j=1,\dots, m,
\end{equation}
with $K$ possibly very large.
The next result shows that this is not possible.
\begin{prop}\label{counterexample}
	Let $\varepsilon >0$, $\alpha,\tilde{c} \geq 0$. Then the following system with $m=2$ and $n=1$
\begin{equation}\label{no-system}
		\left\{\begin{array}{lll}
			u''\pm\varepsilon v'-cu+\alpha v \ge 0 \\
			v''-\tilde c v \ge 0 \quad\quad\quad\quad\quad\quad x\in I_\rho=(0,\rho)
		\end{array} \right.
	\end{equation}
where  $u$,$v$ are scalar functions of $x\in\R$  does not satisfy {\bf wMP}, provided that
\begin{equation}\label{zetacond}
 \zeta(\rho \sqrt c)\sqrt c> \frac{\alpha}{\varepsilon} \quad \mbox{ where } \ \zeta(\tau):= \frac{\cosh  \tau -1}{\sinh \tau-\tau}\;.
 \end{equation}
\end{prop}
\begin{oss}\rm{

	Since $\zeta(0^+)=+\infty$ and $\zeta(+\infty)=1$,
this proposition entails that, for every $\e,K>0$, there exists a system of the type~\eqref{systemcomp}, 
with $B^{(i)}$ satisfying
$|B^{(i)}_{jk}|\leq\e$
and $C=(c_{jk})_{j,k}$ satisfying~\eqref{coop>}, for which {\bf wMP} fails.
Namely, even an 
arbitrary small amount of coupling at the level of first derivatives can prevent the validity of
{\bf wMP} although the zero order matrix is,  so to say, ``very strongly cooperative''.
It also shows that, for any $\varepsilon,c >0$ and $\alpha,\tilde{c} \geq 0$, 
{\bf wMP} fails for~\eqref{no-system} in a small enough interval $I_\rho$.
The fact that {\bf wMP} fails 
when the diagonal zero-order term $c$ is sufficiently large or when 
the size $\rho$ of the interval is sufficiently small can be surprising, if one has in mind the picture for the
scalar equation (where both having a large --negative-- zero-order term and a small domain help
the validity of the maximum principle).\\
This phenomenon could be related to a non-monotonic structure of the system when a first-order coupling is in force.
}
\end{oss}
\begin{oss}\rm{
A few more comments are in order here. We are considering a system with coupled gradients ($\varepsilon >0$). The first part of Proposition \ref{counterexample} says that {\bf wMP} cannot be satisfied  in all bounded domains as soon as $\varepsilon>0$, whatever the amount of cooperativity ($\alpha>0$) is. The second part means that in a fixed interval {\bf wMP} fails for $c$ large enough. In cooperative systems under consideration ($0 \le \alpha \le c$) an excess of coercivity with respect to the coupling ($c$ large compared with $\alpha/\varepsilon$) seems to be responsible for invalidating {\bf wMP}. In particular this is the case in any interval $I_\rho$ when $\alpha=0$.
}
\end{oss}

We exhibit in Proposition~\ref{counterexample-BIS}
below that the same qualitative phenomenon occurs for a larger class of systems.
The proofs of Propositions \ref{counterexample} and \ref{counterexample-BIS} are detailed in Section 2.
}
\begin{prop}\label{counterexample-BIS}
	For every $\e\neq0$, $\tilde c>0$ and $\tilde\e,\alpha,\beta\in\R$,	
	there exists $c>0$ large enough such that the system
	\begin{equation}\label{noMP}
		\left\{\begin{array}{lll}
			u''-\varepsilon v'-cu+\alpha v \ge 0 \\
			v''-\tilde\varepsilon u'-\tilde c v+\beta u \ge 0
		\end{array} \right.
		\quad\text{in }(0,1)
	\end{equation} 
	violates the {\bf wMP}.
	\newline
	\indent
	If, instead, $c>0$ is also fixed in the system~\eqref{noMP}, there exists an interval
	$I\subset(0,1)$ in which such system does not fulfill the {\bf wMP}.	
\end{prop} 

Let us turn now to the positive results.
The study of sufficient conditions for the
validity of the weak Maximum Principle in the form {\bf wMP} in the case where coupling occurs also at the level of first 
or second order derivatives is apparently less explored in literature, 
see however \cite{KreMa1},\cite{KreMa2},\cite{KreMa3}
and also \cite{Miranda},\cite{LiuZhang} for the related issue of maximum norm estimates of the form
$\sup_{x\in\Omega} |u(x)| \leq C\,\sup_{x\in\Omega} |f(x)|$
for solutions $u$ of non-homogeneous systems of equations involving higher order couplings.

The {\bf wMP} property  (\ref{sp}) can be understood in the  framework of the general theory of invariant sets  
introduced by  H.~F.~Weinberger in \cite{Wein} in the context of elliptic and parabolic weakly coupled systems. 
We refer to the recent paper by G.~Kresin and V.~Mazya \cite{KreMa3} where the notion of invariance is thoroughly developed for general systems with couplings at the first and the second order in the case $C\equiv 0$.\\
According to the notion introduced in \cite{Wein},  a set $S\subseteq \R^m$ is invariant for system (\ref {systemcomp})
%\begin{equation} \label{systemhomo}
%\mathrm{Tr}(A_i(x)\nabla^2u_i) + b_i(x)\cdot\nabla u_i+ \sum_{j=1}^m c_{ij}(x)u_j=0,\quad i=1,\dots,m\\
%\end{equation} 
 if  the following property holds
\begin{equation} \label{invariance}
\mathbf{INV}:\qquad u(x)\in S \;\; \hbox{for all}\; x\in \partial \Omega 
\ \  \implies \ \  u(x) \in S \;\; \hbox {for all}\;  x\in \Omega 
\end{equation}
The sign propagation property (\ref{sp}) can then be rephrased as the property of the negative orthant $\R^m_-= \{ u=(u_1,\dots u_m): u_j \leq 0\,\,, \, j=1,\dots,m  \}$
being an invariant set for system (\ref{systemcomp}) of partial differential inequalities.\\
In \cite{Wein} it is proved in particular that {\bf wMP} holds for weakly coupled uniformly elliptic systems
such as
\begin{equation}\label{W}
\mathrm{Tr }(A^j \nabla^2 u_j)+ b^j\cdot \nabla u_j+f(u) =0\;,\; j=1,\dots,m
\end{equation}
under the condition that the vector field $f$ satisfies the property
that for any $p$ belonging to the outward normal cone to $\R^m_-$ at a point $u$ on the boundary of $\R^m_-$ the inequality
\begin{equation} \label{flux}
p\cdot f(u) \leq 0
\end{equation}
holds.
For $f(u)=Cu$, this geometric condition turns out to be the cooperativity property~(\ref{coop}) of matrix $C$. Note also that this condition implies that  $\R^m_-$ is invariant under the flow $du/dt= Cu\,,t >0$.\\

We recall that Proposition~\ref{counterexample} entails that  $\R^m_-$ may fail to be an invariant set
even when the coupling of the first order terms is very small.
As a matter of fact, the first order matrix of system~\eqref{no-system} is $$\left(\begin{array}{cc}0 & -\varepsilon\\0 & 0\end{array}\right)$$
which is not diagonalizable.
This is indeed consistent with results in \cite{KreMa3}. It is in fact shown in that paper, see in particular results in Section 3, that  the sufficient conditions involving the relations between the geometry of a closed convex set $S$ and the matrices
$B^{(i)}$ which imply the invariance of $S$, necessarily require, in the case $S=\R^m_{-}$, the diagonal structure of the first order couplings.

On the account of the example (\ref {no-system}) exhibited in Proposition \ref{counterexample} we are forced to investigate the validity of a weaker form of the sign propagation property or, in other words, to single out an appropriate invariant set for system (\ref{systemcomp}) when first order couplings occur.\\

 It turns out that under some algebraic conditions, including notably the simultaneous diagonalizability of the matrices $B^{(i)}$, a cone propagation type result holds:
 \begin{theorem}\label{main}
Let $\Omega$ be a bounded open subset  of $ \R^n$. Assume that there exists an invertible $m\times m$ matrix $Q$ such that, for all $i=1,\dots,n$,
\begin{equation} \label {diago}
Q^{-1} B^{(i)} Q= \mathrm{Diag}\big(\beta_1^{(i)},\dots,\beta_m^{(i)}\big)\quad \hbox{for some}\;\; \beta^{(i)}_j \in \R, \
\ \ (j=1,\dots,m)
\end{equation} 
\begin{equation}\label{positiv}
Q^{-1}\geq 0
\end{equation}
%\begin{equation}\label{commute}
%AD=DA
%\end{equation}
and, moreover,
\begin{equation}\label {coop2}
Q^{-1}CQ\ \ \hbox{fulfills the cooperativity condition \eqref{coop}}
\end{equation}
\indent
Then the convex cone $ S= \{ u\in \R^m: Q^{-1}u \leq 0\}$ is invariant for system (\ref{systemcomp}).\\
\end{theorem}

%\mathbf{wwMP} \quad D^{-1} u\leq 0  \ \mbox{on} \ \partial \Omega \;\; \mbox{implies}\; D^{-1} u\leq 0 \ \mbox{in}\ \Omega\ \end{equation} 
%for all solutions $u$ of the system of differential inequalities (\ref{system1}).
%Observe that the directional ellipticity condition on $A_i$  is weaker than uniform ellipticity required in previous papers, see \cite{PW},\cite{Sweers},\cite{BMS},\cite{Amann},\cite{BS}.\\

\begin{oss}{ \rm Concerning the linear algebraic conditions of Theorem~\ref{main}, observe first that a matrix $Q$ simultaneously satisfying (1.11) for $i=1,\dots,n$ exists if the $B^{(i)}$'s have a common basis of eigenvectors. This is the case when the matrices $B^{(i)}$ commute each other for all $i=1,\dots,n.$
Observe also that  if $Q$  is an invertible M-matrix, that is $Q=s I- X$ where $X\geq 0$ and $s$ is strictly greater than the spectral radius of $X$, then
$Q$ fulfills condition (\ref{positiv}), see \cite{BP}.\\
Next, it is worth to point out that conditions (\ref{positiv}) and (\ref{coop2}) are compatible.\\ For example, $
 Q=\left(\begin{array}{cc}2 &-1 \\-1 & 2\end{array}\right)
$ is an invertible M-matrix, $
C=\left(\begin{array}{cc}-3&2\\1 & -2\end{array}\right)
$  is cooperative and  
$
 Q^{-1}CQ=\left(\begin{array}{cc}-4&3 \\0 & -1\end{array}\right)
$ is cooperative as well. \\
If no coupling occurs in first derivatives, so that $Q=Q^{-1}=I$, the above result reproduces the one in \cite{CDV3}.
}
\end{oss}

\begin{oss}{ \rm  A related remark is that  permutation matrices satisfies both $Q^{-1}\geq 0$ and $Q\geq 0$, so that in this case the conclusion of Theorem \ref{main} is in fact that the negative orthant  $R^m_-$ is invariant.
 However, it is easy to check that in this situation condition (\ref{diago}) implies that each $B^{(i)}$  is diagonal and the results of \cite{CDV3} apply.\\
A further remark is that one cannot expect in general the invariance of the negative orthant  $R^m_-$. This is indeed coherent with results in \cite{KreMa3};  Lemma 2 there states in fact that the geometric sufficient condition on the matrices $B^{(i)}$ guaranteeing the invariance of $R^m_-$ implies their diagonal structure.   
}
\end{oss}
\begin{oss}{ \rm  Theorem~\ref{main} can in fact be extended (with a completely analogous proof) to a second order matrix operator \begin{equation}\label{A}
Au=\left(\begin{array}{c}\mathrm{Tr }(A \nabla^2 u_1) \\. \\.\\. \\\mathrm{Tr}(A \nabla^2 u_m)\end{array}\right)
\end{equation}
where $A$ is a positive semidefinite matrix such that $A\nu\cdot\nu\ge \lambda>0$ for some direction $\nu\in\R^n$. For some applications of this notion of directional uniform ellipticity condition see \cite{CLN},\cite{CDV1},\cite{CDV2},\cite{VMfe}.
}
\end{oss}
A key role in the proof of this result, which is postponed to the next section, is based on a reduction to a
suitable fully nonlinear scalar differential inequality governed by the elliptic convex Bellman-type 
operator $F$ defined, on scalar functions $\psi:\Omega\to\R$, as
\begin{equation}\label{Bellman}
F[\psi] = \Delta \psi+ \max_{j=1,\dots,m} \sum_{i=1}^n \beta^{(i)}_j \frac{\partial \psi}{\partial x_i} =
\Delta \psi+ \max_{j=1,\dots,m} b^j\cdot\nabla\psi
\end{equation} 
where $\beta^{(i)}_j$ are as in~\eqref{diago}
and $b^j:=(\beta^{(1)}_j,\dots,\beta^{(m)}_j)$.
The main ingredients in the proof are results in \cite{CDV3}, see in particular Theorems 1.1 and 1.3, and the notion of
generalized principal eigenvalue for scalar fully nonlinear degenerate elliptic operators and its relations with  the validity of {\bf wMP}, see \cite{BCPR}.\\

The next example provides a simple illustration of the result of Theorem \ref{main}:
\begin{ex}
{\rm Let $u=(u_1,u_2)$  be a solution of 
	\[
\label{no-system2}
		\begin{cases}
		\displaystyle
			\Delta u_1+6\frac{\partial u_1}{\partial x_1}+\frac{\partial u_2}{\partial x_1}-u_1 \ge 0\\ 
	\displaystyle
			\Delta u_2 - 8\frac{\partial u_1}{\partial x_1}-u_2 \ge 0
		\end{cases}
			\]
in a bounded domain $\Omega\subset\R^n$. In this case
$ B^{(1)}=\left(\begin{array}{cc}6& 1 \\-8 & 0\end{array}\right)$, $B^{(2)}=0$, 
$C=\left(\begin{array}{cc}-1& 0 \\0 & -1\end{array}\right)$
and Theorem  \ref{main} applies with
$$Q=\left(\begin{array}{cc}-1 & 1/2 \\4 & -1\end{array}\right)\quad
 Q^{-1}=\left(\begin{array}{cc}1 & 1/2 \\4 & 1\end{array}\right)
$$ 
yielding  that inequality $u_2\leq \min(-2u_1;-4u_1)$ propagates from $\partial \Omega$ to the whole $\Omega$.
}
\end{ex}
The result of Theorem \ref{main} can be somewhat refined by a suitable weakening of the assumptions there. Firstly, observe that $B^{(i)}$ is not necessarily diagonalizable. A suitable change of basis generally generally leads to an upper triangular matrix, which yields a  real Jordan canonical form of $B^{(i)}$. \\
Suppose that  the $B^{(i)}$\,'s have a common eigenspace of dimension $k \le m$ and consider a basis of $\mathbb R^m$ where the first $k$ vectors are linearly independent (common) eigenvectors of $B^{(i)}$, then  we can find an $m\times m$ invertible real matrix $\hat Q$ that produces a real Jordan canonical form $J^{(i)}=\hat Q^{-1} B^{(i)}\hat Q$, where the $k\times m$ sub-matrix  with the first $k$ rows is made up by a $k\times k$  diagonal block $\Lambda$  plus the $k \times (m-k)$  zero matrix.

In this setting we have the following:

\begin{theorem}\label{main2} Assume in addition to the above that the $k\times m$ sub-matrix containing the first $k$ rows of $\hat C=\hat Q^{-1}C\hat Q$  is made by a cooperative $k\times k$ block  plus the $k \times (m-k)$  zero matrix. Let $p_{ij}$ be the entries of the matrix $\hat P:=\hat Q^{-1}$.\\ If  the $k\times m$ sub-matrix  with the first $k$ rows of $\hat P$ is positive, then the closed convex set
$$
\hat S=\left \{ u\in \R^m:\sum_{j=1}^{m}p_{1j}u_j\leq 0, \dots, \sum_{j=1}^{m}p_{kj}u_j\leq 0 \right\} 
$$
is invariant for the system (1.1).
\end{theorem}
An illustrative example is provided next:
\begin{ex}
{\rm Consider the $2\times 2$ system
	\[
\label{no-system2}
		\begin{cases}
		\displaystyle
			\Delta u_1+\frac{\partial u_2}{\partial x_1}-u_1 \ge 0\\ 
	\displaystyle
			\Delta u_2 +\frac{\partial u_1}{\partial x_1}-u_2 \ge 0
		\end{cases}
			\]
in a domain $\Omega$ of $\R^2$. In this case
$ B^{(1)}=\left(\begin{array}{cc}0& 1 \\1 & 0\end{array}\right)$, $B^{(2)}=0$, 
$C=\left(\begin{array}{cc}-1& 0 \\0 & -1\end{array}\right)$
and Theorem  \ref{main2} applies with
$$\hat Q=\left(\begin{array}{cc}1 & 1 \\1 & -1\end{array}\right)\quad
\hat Q^{-1}=\left(\begin{array}{cc}1/2 & 1/2 \\1/2 & -1/2\end{array}\right)
$$ Since the first row of $\hat P=\hat Q^{-1}$ is nonnegative the above result yields the invariance of  the convex set $S=\{ u=(u^1,u^2): u^1+u^2\leq 0\}$.\\
Note that {\bf wMP}, that is the invariance of $\R^2_{-}$, does not hold true in this example. Indeed, the vector $u=\big((x_1-x_1^2)(x_2-x_2^2)^3,(x_1^2+2x_1-4)(x_2-x_2^2)\big)$
is a solution in the square $\Omega= [0,1] \times[0,1]$ taking non positive values on $\partial \Omega$ with $u_2\leq 0$, $u_1\geq 0$ in $\Omega$.
}
\end{ex}
%%%%%%%%%%%%%%%%%%%%%%%%%%%%%%%%%%%%%%%%%%%%%%%%%%%%%%%%%%%%%%%%%%%%%%%%%%%%%%%%%%%%%%%%%%%%%%%%%%%%%%%%%%%%%%%%%%%%%%%%%%%%%%%%%%%%%%%%%%%%%%%%%%%%%%%%%%%%%%%%%%%%%%%%%%%%%%%%%%%%%%%%%%%%%%%%%%%%%%%%%%%%%%%%%%%%%%%%
\begin{section}{Proofs of the results}
The first part of this section is dedicated to the proofs of Proposition \ref{counterexample} and \ref {counterexample-BIS}.\\
\begin{proof}[Proof of Proposition~\ref{counterexample}]

We restrict to the case $-\varepsilon$, with $\varepsilon>0$. In fact, we can reduce to it by using the change of coordinate $x \to \rho-x$. We also observe that the argument is not affected by $\tilde c\ge 0$, so that we will omit to mention it when discussing on the parameters.

Firstly, we observe that $v(x)=-x$ obviously satisfies the second equation with $v \le 0$ in~$[0,\rho]$.\\
Next, we introduce the sequence of functions
\begin{equation}\label{u}
u_k(x)=-\frac\varepsilon c \left\{ \frac{1-e^{-\sqrt c/k}}{2\sinh(\sqrt c/k)}\,e^{\sqrt c x}+ \frac{e^{\sqrt c/k}-1}{2\sinh(\sqrt c/k)}\,e^{-\sqrt c x}-1-\frac\alpha\varepsilon\left(\frac{\sinh(\sqrt c x)}{k\sinh(\sqrt c/k)}-x\right)  \right\}.
\end{equation}
Then, a direct computation, shows that  for all $k \in \mathbb N$, 
\begin{equation}\label{u-sol}
u_{k}''-\varepsilon v'-cu_k+\alpha v =  0  \ \ \hbox{\rm in} \ I_\rho\equiv (0,\rho)
\end{equation}
and $u_k(0)=0$.

The case $c=0$ is ruled out  either by taking the  limit as $k\to +\infty$ or directly by putting $v=-x$ and $c=0$ in the above equation.

Note also that for $k \to \infty$
\begin{equation}\label{limitk}
\begin{split}
u_{k}(x) \to  u_0(x)&=-\frac\varepsilon c \left\{ \cosh\sqrt c x-1- \frac\alpha\varepsilon \left(\frac{\sinh\sqrt c x}{\sqrt  c}-x\right)\right\}
\end{split}
\end{equation}

Let the parameters $\varepsilon$, $c$, $\alpha$ and $\rho$ be fixed. For  large $k \in \mathbb N$:
\begin{equation}\label{u'}
\begin{split}
u_k'(0)&=-\frac\varepsilon c \left\{ \frac{1-\cosh(\sqrt c/k)}{\sinh(\sqrt c/k)/\sqrt c} -\frac\alpha\varepsilon\left(\frac{\sqrt c/k}{\sinh(\sqrt c/k)}  -1\right) \right\}\\
&=\frac{\varepsilon}{2k}+o(1/k)
\end{split}
\end{equation}
so that $u'_k(0)>0$ for $k$ large enough. 

Since $u_k(0)=0$, we also have $u_k(x)>0$ for some $x \in I_\rho$ for such $k \in \mathbb N$. So {\bf wMP} will be violated if $u_k(\rho)\le 0$.

Next, computing (\ref{limitk}) for $x=\rho$, 
\begin{equation}\label{u0}
\begin{split}
 u_k(\rho)&=-\frac{\varepsilon\sqrt c}{\sinh\sqrt c \rho-\sqrt c \rho} \left\{ \cosh\sqrt c \rho-1- \frac\alpha\varepsilon \left(\frac{\sinh\sqrt c \rho}{\sqrt  c}-\rho\right)\right\}+u_k(\rho)-u_0(\rho)\\
&=-\frac{\varepsilon\sqrt c}{\sinh\sqrt c \rho-\sqrt c \rho} \left\{  \zeta(\sqrt c\rho)\sqrt c- \frac\alpha\varepsilon\right\}+u_k(\rho)-u_0(\rho)
\end{split}
\end{equation}
where 
\begin{equation}\label{zeta}
\begin{split}
 \zeta(\tau)&= \frac{\cosh\tau -1}{\sinh \tau- \tau}.
\end{split}
\end{equation}

Therefore condition $ \zeta(\sqrt c\rho)\sqrt c > \frac \alpha\varepsilon$, see (\ref{zetacond}), yields $u_k(\rho)\le 0$,  for large $k\in \mathbb N$, so that {\bf wMP} is not satisfied. Once established this  fact, we search for condition (\ref{zetacond}) to prove that {\bf wMP} fails.\\
 A straightforward calculation shows that $\zeta(\tau) \to \infty$ as $\tau \to 0^+$ and $\zeta(\tau) \to 1$ as $\tau \to \infty$.
Therefore there exists $\rho_0=\rho_0(c;\frac\alpha\varepsilon)$ such that condition (\ref{zetacond}) holds for $\rho<\rho_0$, and {\bf wMP} is not satisfied, thereby proving that as soon as $\varepsilon >0$ there are intervals $I_\rho$, small enough, where {\bf wMP} fails, whatever $c$ and $\alpha$ are.

On the other hand, let $\rho>0$ be fixed.  The function $\zeta(\sqrt c\rho)\sqrt c$ is increasing with respect to $c$, 
and $\zeta(\sqrt c\rho) \to 1$ as $c \to \infty$, so that
\begin{equation}\label{clarge}
\lim_{c \to \infty}\zeta(\sqrt{c} \rho)\sqrt c=\infty.
\end{equation}
Hence there exists $c_0=c_0(\rho;\frac\alpha\varepsilon)$ such that condition (\ref{zetacond}) holds for $c>c_0$, and {\bf wMP} is not satisfied, thereby proving that as soon as $\varepsilon >0$ then {\bf wMP} fails  in any interval $I_\rho$ and for any $\alpha \ge 0$ when a sufficiently large $c$ is taken.

By the increasing monotonicity of the function $c \to \zeta(\sqrt c\rho)\sqrt c$ we get
\begin{equation*}\label{csmall}
\inf_{c>0}\zeta(\sqrt{c} \rho)\sqrt c=\lim_{c \to 0^+}\zeta(\sqrt{c} \rho)\sqrt c=\frac3\rho\,.
\end{equation*}

It follows that, if $\frac\alpha\varepsilon < \frac3\rho$, then condition (\ref{zetacond}) is satisfied for all $c>0$. This means that in this case we can choose $c_0(\rho;\frac\alpha\varepsilon)=0$.

Finally, recalling that $\zeta(\sqrt c \rho) \to 1$ as $c \to \infty$, then $\zeta(\sqrt c \rho)\sqrt c \cong \sqrt c$ for $c\ge c_1(\rho)$. Hence condition (\ref{zetacond}) is equivalent to 
\begin{equation}\label{clarge}
\sqrt c>\frac\alpha\varepsilon.
\end{equation}
It follows that, if $\frac\alpha\varepsilon>c_1$, then we can choose $c_0(\rho;\frac\alpha\varepsilon)=\left(\frac\alpha\varepsilon\right)^2$.
\end{proof}

\begin{figure}[ht]
\vskip -0,25 cm
\begin{center}
\includegraphics[height=10cm,width=16cm,angle=0]{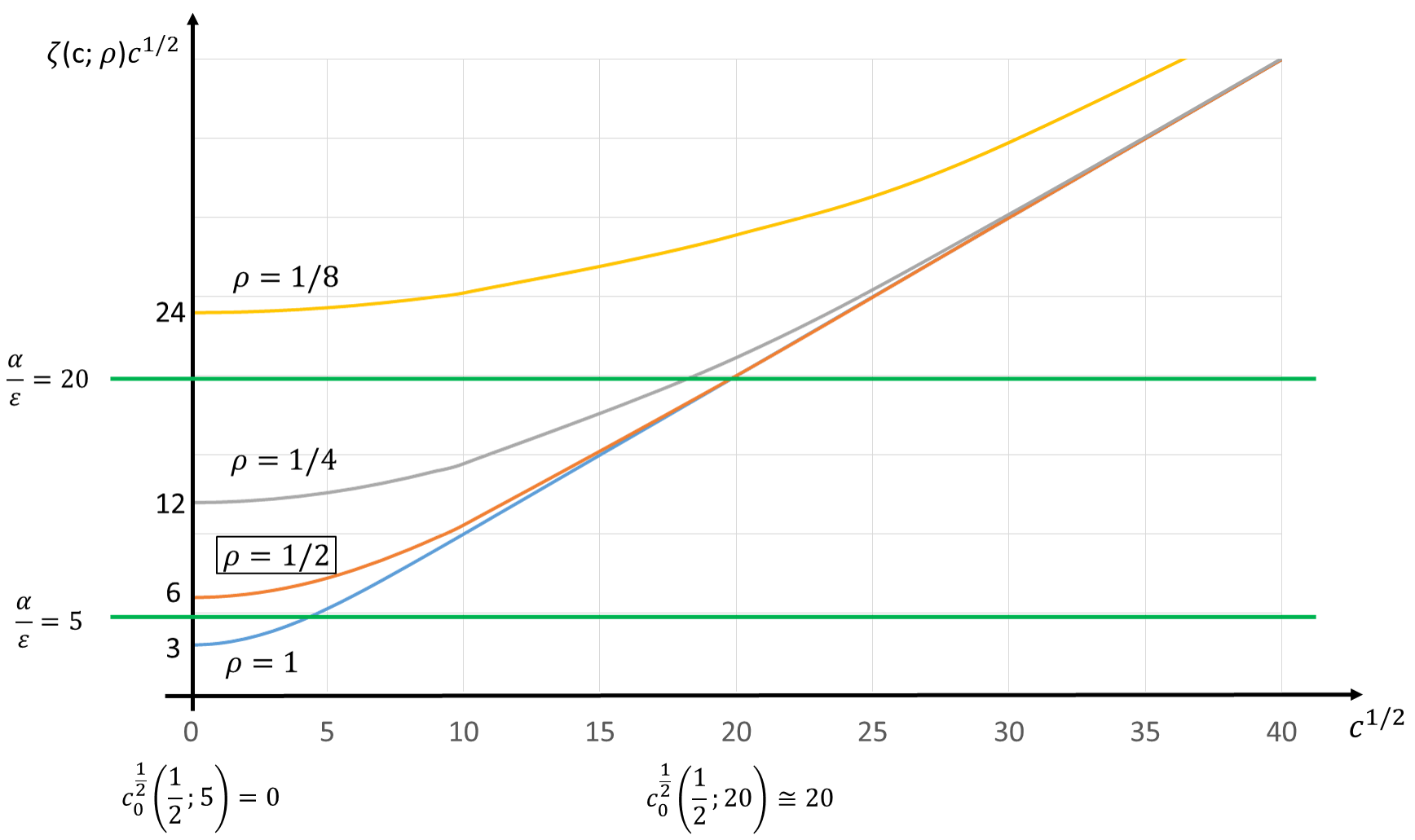}
\vskip 0cm
\caption{The function $\zeta(\sqrt c \rho)\sqrt c$}
\vskip0cm
\label{zeta.fig}
\end{center}
\end{figure}

A picture of the function $c \to \zeta(\sqrt c \rho)\sqrt c$ for different values of $\rho>0$  is shown in Figure \ref{zeta.fig}, where condition (\ref{zetacond}) with the threshold $c_0$ is graphically exhibited on the track $\rho=\frac12$ for different values of $\frac\alpha\varepsilon$.

	\begin{proof}[Proof of Proposition \ref{counterexample-BIS}]
		Up to replacing $u(x),v(x)$ with $u(-x),v(-x)$, it is not restrictive to assume that
		$\e>0$. 
		%	We write $\tilde\e:=-\e>0$.	
		We claim that, for $c$ sufficiently large, there exists a pair $(u,v)$ 
		satisfying~\eqref{noMP} in a strict sense, namely
		\[
		\left\{\begin{array}{lll}
			\inf_{(0,1)}\Big(u''-\varepsilon v'-cu+\alpha v\Big) > 0 \\
			\inf_{(0,1)}\Big(v''-\tilde\varepsilon u'-\tilde c v+\beta u\Big)> 0
		\end{array} \right.
		\] 
		such that $u(0)=0>u(1)$, $v(0)\leq0$, $v(1)\leq0$ and $u'(0)=0$.
		Then, for $\delta>0$ sufficiently small, the pair of functions $(u(x)+\delta x,v(x))$ 
		still satisfies the system~\eqref{noMP} and both functions are $\leq0$ on the boundary of
		$(0,1)$, but $u(x)>0$ for $x>0$ small, hence the {\bf wMP} is violated.
		
		Let us construct the pair of strict subsolutions $(u,v)$. They are defined as follows:
		$$v(x)=x^2-x,\qquad u(x)=\sigma\chi(x),$$
		where  $\chi$ is a smooth, non-increasing function
		satisfying
		$$\chi(0)=\chi'(0)=0,\qquad \chi(x)=-1\ \ \text{for }\,x\geq
		\min\Big\{\,\frac14\,,\,\frac{\e}{4|\alpha|+1}\,\Big\}$$
		and $\sigma$ is a positive constant that will be chosen later.
		We compute, for $x\in(0,1)$,
		$$v''-\tilde\e u'-\tilde c v+\beta u\geq2-\sigma\big(|\beta|+\tilde\e\|\chi'\|_{L^\infty((0,1))}\big),$$
		which is larger than $1$ for $\sigma\leq\sigma_1:= 1/(|\beta|+\tilde\e\|\chi'\|_{L^\infty((0,1))}+1)$.
		Next, for $x\in(0,1)$, we have that
		\begin{align*}
			u''-\e v'-cu+\alpha v &
			\geq \e(1-2x)-|\alpha|x-\sigma|\chi''|-c\sigma\chi.
		\end{align*}
		We estimate the right-hand considering first $0<x\leq\min\{1/4,\e/(4|\alpha|+1)\}$, where we have  
		$$u''-\e v'-cu+\alpha v\geq\frac{\e}4-\sigma\|\chi''\|_{L^\infty((0,1))},$$
		which is larger than $\e/8$ for $\sigma\leq\sigma_2:=\e/(8\|\chi''\|_{L^\infty((0,1))})$.
		While, for $\min(1/4,|\e|/(4|\alpha|+1))\!<\!x\!<\!1$, we see that
		$$u''-\e v'-cu+\alpha v
		\geq -\e-|\alpha|+c\sigma,$$
		which is positive for $c>(\e+|\alpha|)/\sigma$.
		Summing up, taking $\sigma=\min\{\sigma_1,\sigma_2\}$ and then $c>(\e+|\alpha|)/\sigma$,
		we have that $(u,v)$ satisfies~\eqref{noMP} in a strict sense.
		The first statement of the proposition is thereby proved.
		
		Let us turn to the second statement. We have seen above that  {\bf wMP} fails 
		for~\eqref{noMP} provided~$c$ is larger than some $\bar c>0$, and more precisely that 
		it is violated by a pair $(u,v)$ with $v<0$ on~$(0,1)$ and $u>0$ somewhere.
		Consider the pair $(u,v)$ associated with $c=\bar c$
		and let $I$ be a connected component of the set where $u>0$ in $(0,1)$, hence 
		$u=0$ on $\partial I$. For $c\leq\bar c$ there holds in $I$,
		$$u''-\e v'-cu+\alpha v
		\geq u''-\e v'-\bar cu+\alpha v\geq0.$$
		This means that the {\bf wMP} fails in $I$ and then concludes the proof.
	\end{proof}

Let us go now to the proof of Theorem \ref{main}. 
\begin{proof}[Proof of Theorem~\ref{main}]
	 Assume that $u\in [C^2(\Omega)]^m\cap [C^0(\overline\Omega)]^m$ satisfies  (\ref{systemcomp}) and that
	 $u\leq 0$ on $\partial \Omega$. Set 
	 $$\hat B^{(i)}:= Q^{-1} B^{(i)} Q\ \ \ \text{and} \ \ \ \hat C:= Q^{-1} C Q.$$
	 Observe that the change of unknown $u=Qv$ gives, on the account of assumptions (\ref{diago}), (\ref{positiv}), 
	 that $v$ satisfies 
\begin{equation}\label{systemcomp2}
Av+\sum_{i=1}^n \hat B^{(i)} D_iv+\hat C v \geq 0\;\; \hbox{in}\;\; \Omega
\;\; \hbox{and}\;\; v\leq 0 \; \hbox{on}\;\; \partial \Omega
\end{equation}
that is, componentwise,
\begin{equation}\label{system3}
		\left\{\begin{array}{lll}
			\Delta v_1+b^1\cdot \nabla v_1+\hat C_1 v\,\,\,\;\,\, \ge 0 \\
%			\Delta v_2+b^2\cdot \nabla v_2+\hat C_2 v\,\,\,\;\,\, \ge 0\\
			\qquad\qquad\cdots\\
			\Delta v_m+b^m\cdot \nabla v_m+\hat C_m v\ge 0\\
		\end{array} \right.
	\end{equation}
where $b^j=(\beta^{(1)}_j,\dots,\beta^{(m)}_j)$ and $\hat C_j$ is the $j$-th row of $\hat C$,
for $j=1,\dots,m$.
\\
We now employ the argument of the proof of Theorem~1 in~\cite{CDV3} which reduces the above system to a 
scalar inequality governed by the uniformly elliptic (nonlinear) Bellman operator $F$ in (\ref{Bellman}).
By viscosity calculus results based on the cooperativity condition (\ref{coop2}), 
see \cite{CDV3,BS}, since
$v=(v_1,\dots,v_m)$ is a classical solution of (\ref{systemcomp2})
then the scalar function
$$ v^*(x):=  \max_{j=1,\dots,m} (v_j)^+(x)\,,$$ where ``$^+$'' denotes the positive part, is a continuous weak solution in the viscosity sense, 
see \cite{CIL},~of 
\begin{equation}\label{Bellman2}
 F[v^*] \geq 0 \;\; \hbox{in}\;\; \Omega
\;\; \hbox{and}\;\; v^*=0 \; \hbox{on}\;\; \partial \Omega\,.
 \end{equation}

Suppose indeed that a smooth function $\varphi$ touches from above $v^*$ at some point in~$\Omega$.
If at that point $v^*=0$ then clearly $F[\varphi]\geq0$ there. Otherwise
$\varphi$ touches from above the component $v_j$ realizing the positive maximum $v^*$ at that 
point and thus there holds
$$\Delta\varphi+b^j\cdot\nabla\varphi+\hat C_j v\geq0.$$
But then recalling that $\hat C_j$ fulfills the cooperativity condition (\ref{coop}),
one infers that
$$\hat C_j v\leq v_j\sum_k\hat C_{jk}\leq0,$$
whence again $F[\varphi]\geq0$.

In order to apply the general result of~\cite{BCPR} we need to show that the generalized principal eigenvalue, see \cite{BCPR},
of $F$ is positive, which amounts to finding a strict supersolution which is strictly positive in~$\overline\Omega$.
The latter is simply provided by $\psi(x)=\psi(x_1,\dots,x_m) = 1 - \delta e^{\gamma x_1}$.
Indeed, this function satisfies
$$F[\psi]=-\delta\gamma e^{\gamma x_1}\Big(\gamma+\min_{j=1,\dots,m} \beta_j^{(1)}\Big),$$
which is strictly negative in $\R^n$ provided $\gamma>|\min_{j=1,\dots,m} \beta_j^{(1)}|$.
We then choose $\delta$ small enough, depending on $\gamma$ and $\Omega$, so that $\psi>0$ in $\overline\Omega$. 
Summing up, $\psi$ is positive in $\overline\Omega$ 
and satisfies there $F[\psi]<0$, hence also $F[\psi]+\lambda\psi<0$ for $\lambda>0$ suitably small. 
This implies 
%by results in \cite{CDV3}
 that the numerical index $\mu_1(F,\Omega)$ defined by
\begin{equation}\label{mu1}
 \mu_1(F,\Omega)= \sup\{\lambda\in\R: \psi\in C(\overline\Omega),\  \psi>0,\ F[\psi ]+\lambda \psi \leq 0 \;\mbox{in}\; \Omega\}
\end{equation}
is strictly positive. Therefore, according to \cite{BCPR}, the weak Maximum Principle for the scalar problem (\ref{Bellman2}) holds, that is $v^*\leq 0$  in $\Omega$.

% At this point Theorem 1.3 in \cite{CDV3} applies giving $(v_j)^+\leq 0$ in $\Omega$ so that 
This means that $Q^{-1}u=v\leq 0$ in $\Omega$ and the proof is complete.
\end{proof}

We conclude the section with the proof of Theorem \ref{main2}
\begin{proof}[Proof of Theorem~\ref{main2}]
Following the same lines of the proof of Theorem  \ref{main}, we set $u=\hat Qv$. When multiplying by $\hat P=\hat Q^{-1}$, this time  we keep, by assumption, the positivity for the first $k$ equations, which again by the assupmtions made are decoupled in the gradient variables. So, letting $\tilde C$ be the diagonal part of $\hat C$ and $\tilde v=(v_1,\dots,v_k)$, we get
\begin{equation}\label{systemk}
		\left\{\begin{array}{lll}
			\Delta v_1+b^1\cdot \nabla v_k+\tilde C_1 \tilde v\ge 0 \\
%			\Delta v_2+b^2\cdot \nabla v_2+\hat C_2 v\,\,\,\;\,\, \ge 0\\
			\qquad\qquad\cdots\\
			\Delta v_k+b^m\cdot \nabla v_k+\tilde C_k \tilde v\ge 0\\
		\end{array} \right.
	\end{equation}
where  $\tilde C_j$ is the $j$-th row of $\hat C$.
The conclusion follows as in the proof of Theorem 2 with $k$ instead of $m$.
\end{proof}

\end{section}

\end{document}